\documentclass[11pt]{article}

\usepackage{amsmath,amssymb}
\usepackage{graphicx}

\newtheorem{theorem}{Theorem}
\newtheorem{ass}{Theorem}[section]
\newtheorem{prop}[ass]{Proposition}
\newtheorem{lemma}[ass]{Lemma}
\newtheorem{claim}[ass]{Claim}
\newtheorem{corollary}[ass]{Corollary}

\newtheorem{definition}[ass]{Definition}

\newcommand{\qed}{\hspace*{\fill} \rule{7pt}{7pt}}
\newcommand{\Proof}{\noindent{\bf Proof.}\ \ }

\begin{document}

\title{On the maximum number of spanning copies of an orientation in a tournament}

\author{
Raphael Yuster
\thanks{Department of Mathematics, University of Haifa, Haifa
31905, Israel. Email: raphy@math.haifa.ac.il}
}

\date{}

\maketitle

\setcounter{page}{1}

\begin{abstract}

For an orientation $H$ with $n$ vertices, let $T(H)$ denote the maximum possible number of labeled copies of $H$ in an $n$-vertex tournament.
It is easily seen that $T(H) \ge n!/2^{e(H)}$ as the latter is the expected number of such copies in a random tournament.
For $n$ odd, let $R(H)$ denote the maximum possible number of labeled copies of $H$ in an $n$-vertex regular tournament.
Adler et al. proved that, in fact, for $H=C_n$ the directed Hamilton cycle, $T(C_n) \ge (e-o(1))n!/2^{n}$ and it was observed by Alon that already
$R(C_n) \ge (e-o(1))n!/2^{n}$. Similar results hold for the directed Hamilton path $P_n$.
In other words, for the Hamilton path and cycle, the lower bound derived from the expectation argument can be improved by a constant factor.
In this paper we significantly extend these results and prove that they hold for a larger family of orientations $H$ which includes all bounded degree Eulerian orientations
and all bounded degree balanced orientations, as well as many others. One corollary of our method is that for any $k$-regular orientation $H$ with $n$ vertices,
$T(H) \ge (e^k-o(1))n!/2^{e(H)}$ and in fact, for $n$ odd, $R(H) \ge (e^k-o(1))n!/2^{e(H)}$.
\end{abstract}

\section{Introduction}

All graphs and digraph considered in this paper are finite and simple. Our graph-theoretic notation is standard and follows that of \cite{bollobas-1978}.
An {\em orientation} is a digraph without directed cycles of length $2$.
Stated otherwise, an orientation $H$ of an undirected graph $\hat{H}$ is obtained by
orienting each edge of $\hat{H}$.  We say that $\hat{H}$ is the {\em underlying graph} of $H$. 
We use the notation $(u,v)$ for an edge  directed from $u$ to $v$.
For a vertex $v$ of an orientation its {\em out-degree}, denoted $d^+(v)$, is the number of edges directed from $v$
while its {\em in-degree}, denoted $d^-(v)$, is the number of edges directed to $v$.
An orientation $H$ is {\em $k$-regular} if $d^+(v)=d^-(v)=k$ for all vertices $v \in V(H)$.
It is Eulerian if $d^+(v)=d^-(v)$ for all vertices $v \in V(H)$, and $\hat{H}$ is connected
(if connectivity of $\hat{H}$ is not assumed  we say that $H$ is {\em even}). 
An orientation is {\em balanced} if $|d^+(v)-d^-(v)| \le 1$ for all vertices $v \in V(H)$.

A {\em tournament} is an orientation of a complete graph. 
A tournament is {\em transitive} if it contains no directed cycle, or equivalently, if it is possible to order its vertices
so that all edges are oriented from left to right. We use $T_n$ to denote the (unique) $n$-vertex transitive tournament.
A tournament with $n$ vertices is {\em regular} if $d^+(v)=d^-(v)=(n-1)/2$ for all vertices $v \in V(H)$ (so here $n$ is odd).

For an orientation $H$ and a tournament $G$, both with the same number of vertices,  let
$G(H)$ denote the number of labeled copies of $H$ in $G$.
Notice that each such labeled copy is a labeled spanning subgraph of $G$ and notice also that the
number of unlabeled copies of $H$ in $G$ is $G(H)/aut(H)$ where $aut(H)$ is the number of automorphisms of $H$.
Our objects of interest in this paper are $T(H)$, which denotes the maximum of $G(H)$ ranging over all tournaments
with the same number of vertices of $H$, and $R(H)$, which denotes the maximum of $G(H)$ ranging over all regular tournaments
with the same number of vertices of $H$ (so for $R(H)$ we assume that the number of vertices is odd).

Let $C_n$ denote the directed cycle on $n$ vertices and let $P_n$ denote the directed path on $n$ vertices.
The problems of determining $T(C_n),R(C_n),T(P_n),R(P_n)$ are well-studied and fairly well-understood.
The first lower bound for $T(P_n)$ was proved by Szele \cite{szele-1943} who showed that $T(P_n) \ge n!/2^{n-1}$.
It seems that his argument is the first application of the probabilistic method in combinatorics. Indeed,
observe that $n!/2^{n-1}$ is the expected number of copies of $P_n$ in a random tournament.
Thus, Szele's argument actually proves that $T(H) \ge n!/2^{e(H)}$ for any orientation $H$ with $n$ vertices (we call this the ``expectation lower bound'')
and in particular, $T(C_n) \ge n!/2^{n}$, or, equivalently, since $aut(C_n)=n$, there are tournaments which contain at least $(n-1)!/2^n$ Hamilton cycles.
Szele also proved the upper bound $T(P_n) = O(n!/2^{0.75n})$. Many years later, this upper bound was significantly improved by Alon \cite{alon-1990}
who showed that $T(P_n) = O(n^{3/2} n!/2^{n-1})$. In other words, this upper bound is only larger than the expectation lower bound by a small polynomial factor.
Later, the polynomial was further improved to slightly less than $n^{5/4}$ by Friedgut and Kahn \cite{FK-2005}.
Both of these results apply also to $T(C_n)$.

The first, perhaps surprising improvement over the expectation lower bound, was obtained by Adler, Alon, and Ross \cite{AAR-2001}. They have shown, that, in fact
$T(P_n) \ge (e-o(1))n!/2^{n-1}$ and similarly, $T(C_n) \ge (e-o(1))n!/2^n$. It
was even more surprising that later Wormald \cite{wormald-p} proved that the factor $e$ can be replaced with a constant larger than $2.855$.
His clever probabilistic construction involves a computer assisted verification.
Thomassen (see \cite{bondy-1995}) asked in 1990 whether $R(C_n) \ge n!/2^n$, that is, whether the expectation lower
bound also applies to the restricted case of regular tournaments. Alon \cite{alon-p} observed that
the proof from \cite{AAR-2001} can, in fact, be adjusted to show that it holds also for
regular tournaments, namely, $R(P_n) \ge (e-o(1))n!/2^{n-1}$ and $R(C_n) \ge (e-o(1))n!/2^n$.
In fact, Wormald's proof with the improved constant $2.855$ also applies to the regular setting.
It was suggested in \cite{AAR-2001,FK-2005} that $T(P_n)$ and $R(P_n)$ may in fact be larger than the expectation lower bound only by a constant factor.
Wormald \cite{wormald-p} conjectures that this factor is not larger than $2.856$.
It is appropriate to mention here a related result of Cuckler \cite{cuckler-2007}, who proved that {\em every} regular tournament on $n$ vertices contains at least $n!/(2 + o(1))^n$ Hamilton cycles, thereby resolving a conjecture of Thomassen \cite{thomassen-1985}. Recently, Ferber, Krivelevich and Sudakov \cite{FKS-p} significantly generalized Cuckler's result showing that
an analogous result holds when the host graph is not a tournament but any regular orientation with degree at least $cn$ where $c > 3/8$,
namely, there are $((c+o(1))n/e)^n$ Hamilton cycles. Finally, we refer the interested reader to an excellent survey of Kuhn and Osthus \cite{KO-2012} on Hamilton cycles in directed graphs.

The moral of the results of \cite{AAR-2001} and \cite{wormald-p} is that $T(H)$ and even $R(H)$ ``beat'' the expectation lower bound by a constant factor
for $H=C_n$ and $H=P_n$. Hence, it is of interest to determine which other orientations $H$ have the property that $T(H)$ and $R(H)$ are larger than
$n!/2^{e(H)}$ by a constant factor. Clearly, not every $H$ has this property. For let $M_t$ be the orientation of a perfect matching of $t$ edges.
It is straightforward to verify that $G(M_{n/2})=n!/2^{n/2}$ for every tournament $G$ with $n$ vertices. It seems that the regularity of $C_n$ and the almost regularity of $P_n$
is what makes them prone to appear more than the expected number of times in some carefully designed (regular) tournaments. Our main result in fact proves this fact
in a very general setting which we next define. As it turns out, this setting is general enough to include all Eulerian orientations of bounded maximum degree as well as all ``nontrivial'' balanced
orientations of bounded maximum degree.

For an orientation $H$, let $plus(H)$ be the number of directed paths of length $2$ in $H$ and let $minus(H)$ be the number of anti-directed paths of length $2$. Namely,
$$
plus(H) = \sum_{v \in V(H)} d^+(v)d^{-}(v)~~,~~minus(H) = \sum_{v \in V(H)} \binom{d^+(v)}{2}+\binom{d^{-}(v)}{2}\;.
$$
For an orientation $H$, recall that $\hat{H}$ denotes the underlying graph of $H$. We define the maximum degree of $H$ as $\Delta(H)=\Delta(\hat{H})$.
For $\epsilon > 0$ and a positive integer $k$ we say that $H$ is an {\em $(\epsilon,k)$-consistent orientation} if $H$ has $n$ vertices,
$\delta(H) \le k$ and $plus(H)-minus(H) \ge \epsilon n$. The following is our main result.
\begin{theorem}\label{t:main}
Let $\epsilon > 0$ and let $k$ be a positive integer. There exists $N=N(\epsilon,k)$ and $\delta=\delta(\epsilon,k)$ such that
for every odd $n > N$ and every $(\epsilon,k)$-consistent orientation $H$ with $n$ vertices, $R(H) \ge (1+\delta)n!/2^{e(H)}$.
\end{theorem}
When $n$ is odd, we trivially have $R(H) \le T(H)$. When $n$ is even, $R(H)$ is not defined, but we can define it as the maximum number of labeled copies of $H$
in an $n$-vertex balanced tournament (i.e. tournaments with minimum in-degree $n/2-1$). With this definition, $R(H) \le T(H)$ regardless of the parity of $n$.
The following can be proved as a corollary of our main result.
\begin{corollary}\label{c:main}
Let $\epsilon > 0$ and let $k$ be a positive integer. There exists $N=N(\epsilon,k)$ and $\delta=\delta(\epsilon,k)$ such that
for every $n > N$ and every $(\epsilon,k)$-consistent orientation $H$ with $n$ vertices, $T(H) \ge R(H) \ge (1+\delta)n!/2^{e(H)}$.
\end{corollary}

As mentioned earlier, the family of $(\epsilon,k)$-consistent orientations is quite large, and contains, for example, all Eulerian orientations of maximum degree $k$ using $\epsilon=1$
as well as all {\em balanced orientations} with maximum degree $k$ and not too many isolated vertices, using a sufficiently small $\epsilon$, as the following proposition shows.
\begin{prop}\label{p:euler}
Let $H$ be a balanced orientation with maximum degree $k$ and with at least $(0.5+\epsilon)n$ edges. Then $H$ is a $(2\epsilon/k,k)$-consistent orientation.
Let $H$ be an Eulerian orientation with maximum degree $k$. Then $H$ is a $(1,k)$-consistent orientation.
\end{prop}
We note that we need a balanced orientation to have more than $0.5$ edges in order to apply Theorem \ref{t:main} or Corollary \ref{c:main} due to the above example which shows
that $M_{n/2}$, which is a balanced orientation with $n/2$ edges, has $T(M_{n/2})=n!/2^{n/2}$.

Interestingly, our proof can be optimized to yield relatively large values for $\delta$ when $H$ is an Eulerian orientation.
This is best demonstrated in the case of $k$-regular orientations for which we have the following theorem.
\begin{theorem}\label{t:k-reg}
Let $H$ be a $k$-regular orientation with $n$ vertices. Then $R(H) \ge (e^k-o(1))n!/2^{nk}$.
\end{theorem}
Notice that the case $k=1$ in Theorem \ref{t:k-reg} implies the theorem in \cite{AAR-2001}.
Let us outline in very high level the similarities and differences between the proof in \cite{AAR-2001} and the proof of Theorem \ref{t:main}.
The idea in \cite{AAR-2001} is to fix some maximal packing of $K_n$ with edge-disjoint triangles (packing almost all edges), and then examine the probability space
of random tournaments obtained by orienting each triangle as a directed triangle in one of the two possible directions. A Hamilton cycle (corresponding to
a permutation) either intersects each triangle in one or two edges, and it is possible to estimate rather precisely, using Brun's Sieve, the probability that all its $n$ edges
agree with the orientation of the random tournament. In our proof of Theorem \ref{t:main} we also fix some maximal packing of $K_n$ into edge-disjoint cliques.
However, we cannot take these cliques to be just triangles for the sheer reason that $\hat{H}$ many have many triangles and there is very high probability that
in a random permutation, few of these triangles in $\hat{H}$ will actually identify with a triangle of the packing. As some (or even all) of these triangles might be $T_3$'s
in $H$, there is very high chance that the permutation will be inconsistent with {\em every} orientation of the random tournament.
In fact, we will need to take our cliques in the packing to be $K_t$'s where $t$ is very large, but still constant. But now, there are many more possibilities on how a random permutation
might intersect each $K_t$ and using Brun's sieve to estimate the probability that all edges of $H$ are consistent with the random tournament becomes very difficult.
Instead we use a different method based on the inequality of arithmetic and geometric means which reduces the computation to a feasible one.

In the following section we prove some basic lemmas which are needed for the proof of the main result, Theorem \ref{t:main}. Section 3 then contains the proof of Theorem \ref{t:main}.
In Section 4 we consider balanced and Eulerian orientations and prove Proposition \ref{p:euler} and Theorem \ref{t:k-reg}. The final section contains some concluding remarks
and in particular, it is shown how to adjust Theorem \ref{t:main} in order to obtain Corollary \ref{c:main}.

\section{Lemmata}

Let $G$ be a graph on vertex set $[n]$. A {\em decomposition} of $G$ is a family ${\cal D}$ of pairwise edge-disjoint labeled subgraphs of $G$ such that
each edge of $G$ is covered by an element of ${\cal D}$. If all elements of ${\cal D}$ are isomorphic to the same graph $H$ then ${\cal D}$ is called
an $H$-decomposition. The following is a special case of a seminal result of Wilson \cite{wilson-1975}.

\begin{lemma}\label{l:wilson}
Let $t \ge 3$ be an integer. There exists $N=N_{\ref{l:wilson}}(t)$ such that for all $n > N$, if $n \equiv 1,t \bmod t(t-1)$, then $K_n$ has a $K_t$-decomposition.
\end{lemma}
We note that the requirement that $n \equiv 1,t \bmod t(t-1)$ is a trivial necessary condition for the existence of a $K_t$-decomposition as clearly we must have that $t-1$
divides $n-1$ and that $\binom{t}{2}$ divides $\binom{n}{2}$. When $n$ is not of this form we cannot have a $K_t$-decomposition but we can still have a very large set of edge-disjoint $K_t$'s.
In fact, a result of Gustavsson \cite{gustavsson-p} states that if $G$ is a graph on $n$ vertices such that $e(G)$ is divisible by $\binom{t}{2}$ and that each vertex degree is divisible by
$t-1$, then $G$ will have a $K_t$-decomposition if we assume that the minimum degree of $G$ is sufficiently large. Gustavsson's result also follows as a special case of a seminal result of
Keevash on the existence of designs \cite{keevash-p}. We need here the following even more special case.

\begin{lemma}\label{l:gust}
Let $t \ge 3$ be an integer and let $d \ge 1$. There exists $N=N_{\ref{l:gust}}(t,d)$ such that for all $n > N$ he following holds. Suppose $G$ is a graph with $n$ vertices,
minimum degree at least $n-d$ and suppose furthermore that $e(G)$ is divisible by $\binom{t}{2}$ and that each vertex degree is divisible by $t-1$.
Then $G$ has a $K_t$-decomposition.
\end{lemma}
We note that the best known constants for $N=N_{\ref{l:gust}}(t,d)$ are rather huge. In fact, in our proof of Theorem \ref{t:main} we could have settled for a weaker result
than Lemma \ref{l:gust} which only guarantees that $K_n$ has a packing into edge-disjoint copies of $K_t$ which covers all but $o(n^2)$ edges of $K_n$.
This weaker statement, which actually follows as a special case from a hypergraph packing result of R\"odl \cite{rodl-1985}, works with a much smaller $N$.
However, using this weaker variant would have made the proof of Theorem \ref{t:main} considerably more involved, due in part to the fact that we would still need to show that the leftover
uncovered edges can be decomposed into constant size graphs.

\begin{lemma}\label{l:dec}
Let $t \ge 3$ be an odd integer. There exists $N=N_{\ref{l:dec}}(t)$ such that for all odd $n > N$ the following holds.
There exists a graph $B$ on $n$ vertices such that $B$ has maximum degree at most $3t-5$ and $B$ can be edge-decomposed into at most $n(t-3)/6$ triangles, at most $t-3$ cycles of length $4$
and at most $t-1$ copies of $K_{2t-1}$. Furthermore, $K_n$ can be edge-decomposed into copies of $K_t$ and at most one copy of $B$.
\end{lemma}
\Proof
Let $N_{\ref{l:dec}}(t) = N_{\ref{l:gust}}(t,3t-4)$.
Suppose that $t \ge 3$ is odd and $n > N_{\ref{l:dec}}(t)$ is also odd.
Let $q \equiv n \bmod t-1$ where $q \in \{1,3,5,\ldots,t-2\}$.

Our first step is performed in $(q-1)/2$ iterations (if $q=1$ the first step in ignored).
In the first iteration, we remove from $K_n$ a set of pairwise vertex-disjoint triangles and at most two cycles of length four that cover each vertex precisely once.
This is trivially doable as each number larger than $5$ can be written as $3x+4y$ where $x$ and $y$ are nonnegative integers and $y \le 2$.
After the first iteration the remaining graph $G_1$ has $n$ vertices and each degree is $n-3$.
In the second iteration (if there is one), we do the same, removing from $G_1$ a set of pairwise vertex-disjoint triangles and at most two cycles of length four that cover each vertex precisely once.
This is again doable (less trivially now) using, say, a theorem of Corr\'adi and Hajnal \cite{CH-1963} (with lots of room to spare as this theorem actually enable us do
$\Theta(n)$ iterations, while we only need constantly many iterations).
After the last iteration (iteration $(q-1)/2$) the remaining graph $G'$ has $n$ vertices and each degree is $n-q$. The graph of the deleted edges is, by construction,
a regular graph $B'$ with $n$ vertices and degree $q-1$, edge-decomposed into at most $n(q-1)/6 \le n(t-3)/6$ triangles and at most $q-1\le t-3$ cycles of length $4$.

Notice now that each vertex degree of $G'$ is divisible by $t-1$. However, the number of edges of $G'$ is not necessarily divisible by $t(t-1)/2$.
Nevertheless, $e(G')$ is divisible by $(t-1)/2$ since
Indeed
$$
e(G') = \binom{n}{2} - e(B') = \binom{n}{2} - \frac{n(q-1)}{2} = \frac{n(n-q)}{2} \equiv 0 \bmod \frac{t-1}{2}\;.
$$
The second step further removes several pairwise vertex-disjoint copies of $K_{2t-1}$ until the number of edges of the remaining graph is divisible by $t(t-1)/2$.
To show this is possible we first observe that since $t$ is odd,
$$
gcd\left(e(K_{2t-1}),\frac{t(t-1)}{2}\right) = gcd\left((2t-1)(t-1),\frac{t(t-1)}{2}\right) = \frac{t-1}{2}\;.
$$
So indeed after removing at most $t-1$ copies of $K_{2t-1}$  the number of edges becomes divisible by $t(t-1)/2$.
We can remove these many copies vertex-disjointly since the degree of $G'$ is $n-q$, so Tur\'an's Theorem \cite{turan-1941} lets us do this deletion greedily (with lots of room to spare).
Finally, we arrive at the remaining graph $G$ which has $n$ vertices, its number of edges is divisible by $t(t-1)/2$, each vertex degree of $G$ is either
$n-q$ or $n-q-(2t-2)$ so each vertex degree is divisible by $t-1$. Furthermore, the minimum degree of $G$ is at least $n-q-2t+2 \ge n-3t+4$.
Thus, by Lemma \ref{l:gust}, $G$ has a $K_t$-decomposition. The set of removed edges is a graph $B$ that is edge-decomposed into at most $n(t-3)/6$ triangles, at most $t-3$
cycles of length $4$ and at most $t-1$ copies of $K_{2t-1}$. Furthermore, the maximum degree of $B$ is at most $(n-1)-\delta(G)=3t-5$.
\qed

\begin{definition}\label{d:adjust}(adjusted $t$-decomposition)
For an odd integer $t \ge 3$ and an odd integer $n$, we say that a decomposition ${\cal D}$ of $K_n$ is an {\em adjusted $t$-decomposition}
if the set of edges of all elements of ${\cal D}$ that are not isomorphic to $K_t$ form a graph $B$ on $n$ vertices with maximum degree at most $3t-5$.
Furthermore, the elements of ${\cal D}$ that decompose $B$ are either at most $n(t-3)/6$ triangles, at most $t-3$ cycles of length $4$ and at most $t-1$ copies of $K_{2t-1}$.
\end{definition}
Lemma \ref{l:dec} proves that for all odd $t \ge 3$ and all sufficiently large odd $n$, an adjusted $t$-decomposition exists.

\begin{definition}\label{d:atypical}(atypical copy, typical copy)
Let $t \ge 3$ be an odd integer an let ${\cal D}$ be an adjusted $t$-decomposition of $K_n$ ($n$ odd).
Let ${\cal D}'$ denote the subset of elements of ${\cal D}$ that are not copies of $K_t$.
Let $H$ be an undirected graph with $n$ vertices. We say that a labeled copy of $H$ in $K_n$ is {\em atypical}
if at least one of the following holds:\\
(i) some element of ${\cal D}'$ contains two (or more) edges of $H$,\\
(ii) some $K_t$ of ${\cal D}$ contains three (or more) edges of $H$ where the three edges do not induce a triangle in $H$.\\
A labeled copy of $H$ which is not atypical is {\em typical}.
\end{definition}
Observe that if $H$ is dense, then many of the labeled copies are atypical. For example, if $H=K_n$, then all $n!$ labeled copies are atypical (unless $t=3$ and
${\cal D}$ is a Steiner triple system). The following lemma proves that if $H$ is sparse, then most labeled copies of $H$ will be typical.

\begin{lemma}\label{l:typical}
Let $t \ge 3$ be an odd integer an let ${\cal D}$ be an adjusted $t$-decomposition of $K_n$ ($n$ odd).
Let $H$ be an undirected graph with $n$ vertices and maximum degree $d$. Then there are at most $11d^3t^4(n-1)!$ atypical labeled copies of $H$ in $K_n$.
\end{lemma}
\Proof
It is equivalent, although more convenient, to upper bound the probability that a random labeled copy of $H$ is atypical.
So for $\sigma \in S_n$ drawn at random, we let $H_\sigma$ correspond to the random labeled copy of $H$.
Let us first estimate the probability that condition (i) in Definition \ref{d:atypical} holds.
Let ${\cal D}'$ denote the subset of elements of ${\cal D}$ that are not copies of $K_t$.
Recall that the elements of ${\cal D}'$ are either triangles, $4$-cycles or copies of $K_{2t-1}$.
Fix some triangle $\{x,y,z\}$ of ${\cal D'}$. What is the probability that $H_\sigma$ contains the edge $xy$ and also the edge $yz$?
Given that the vertex of $H$ mapped to $y$ is, say, $u$, this occurs only if two neighbors of $u$ in $H$ map to $x$ and $z$. Since the maximum degree of $H$
is $d$, this occurs with probability at most $d(d-1)/((n-1)(n-2))$. Hence, the probability that $\{x,y,z\}$ contains two (or three) edges of $H_\sigma$
is at most $3d(d-1)/((n-1)(n-2))$. Since ${\cal D}$  is an adjusted $t$-decomposition, ${\cal D}'$ has at most $n(t-3)/6$ triangles.
So the probability that some triangle of ${\cal D}'$ contains two or more edges of $H_\sigma$ is at most
\begin{equation}\label{e:pr-c3}
\frac{n(t-3)}{6} \cdot \frac{3d(d-1)}{(n-1)(n-2)}\;.
\end{equation}
The overall number of edges in the $4$-cycles of ${\cal D}'$ is at most $4(t-3)$ and the overall number of edges in the $K_{2t-1}$ copies of ${\cal D}'$
is at most $(2t-1)(t-1)^2$. The probability of any one of these $4(t-3)+(2t-1)(t-1)^2$ edges to belong to $H_\sigma$ is
$e(H)/\binom{n}{2}$. Since $e(H) \le nd/2$, the probability that at least one of the edges of a $4$-cycle of ${\cal D}'$ or a $K_{2t-1}$ element of ${\cal D}'$
will be in $H_\sigma$ is at most
\begin{equation}\label{e:pr-c4+}
\left(4(t-3)+(2t-1)(t-1)^2\right) \frac{d}{n-1}\;.
\end{equation}
Summing (\ref{e:pr-c3}) and (\ref{e:pr-c4+}) we obtain that the probability that some element of ${\cal D}'$ contains two (or more) edges of $H_\sigma$ is
at most
\begin{equation}\label{e:bad}
\frac{n(t-3)}{6} \cdot \frac{3d(d-1)}{(n-1)(n-2)} + \left(4(t-3)+(2t-1)(t-1)^2\right) \frac{d}{n-1} < \frac{8t^3d^2}{n}\;.
\end{equation}

We next compute the probability that some $K_t$ copy of ${\cal D}$ contains three or more edges of $H_\sigma$, where these three edges do not induce a triangle.
Consider three edges $e,f,g$ of $H$ which do not form a triangle. Then they span either $4$, $5$ or $6$ vertices of $H$.
Consider first the case that they span $6$ vertices, so $e,f,g$ form a matching. Given that $e$ is placed in some $K_t$ element of ${\cal D}$, the probability that
also $f$ and $g$ fall in this element is precisely the probability that their four endpoints are mapped to vertices of this element which is
$$
\frac{(t-2)(t-3)(t-4)(t-5)}{(n-2)(n-3)(n-4)(n-5)} < \frac{t^4}{n^4}\;.
$$
Since the number of choices for $e,f,g$ in $H$ is less than $e(H)^3 < d^3 n^3$, the probability that any three such edges fall in the same $K_t$ copy of ${\cal D}$ is at most
\begin{equation}\label{e:6v}
d^3n^3 \cdot \frac{t^4}{n^4}  \le \frac{d^3t^4}{n}\;.
\end{equation}
Consider next the case that $e,f,g$ span $5$ vertices. So, we must have that, say $e$ is disjoint from $f$ and $g$, and $f,g$ share an endpoint.
Given that $e$ is in some $K_t$ element of ${\cal D}$, the probability that
also $f$ and $g$ fall in this element is precisely the probability that their three endpoints are mapped to vertices of this element which is
$$
\frac{(t-2)(t-3)(t-4)}{(n-2)(n-3)(n-4)} < \frac{t^3}{n^3}\;.
$$
Since the number of choices for $e$ is $e(H)$ and the number of choices for $f,g$ sharing an endpoint is less than $nd^2$, the probability that any three such edges fall in the same $K_t$ copy of ${\cal D}$ is at most
\begin{equation}\label{e:5v}
d^3n^2 \cdot \frac{t^3}{n^3}  \le \frac{d^3t^3}{n}\;.
\end{equation}
Consider last the case that $e,f,g$ span $4$ vertices (so either a path on three edges or a star on three edges). 
Given that $e$ is in some $K_t$ element of ${\cal D}$, the probability that
also $f$ and $g$ fall in this element is precisely the probability that the two vertices of $f$ and $g$ that are not endpoints of $e$ are mapped to vertices of this element which is
$$
\frac{(t-2)(t-3)}{(n-2)(n-3)} < \frac{t^2}{n^2}\;.
$$
Since the number of choices for stars with three edges in $H$ is less than $nd^3/2$ and the number of choice of paths on three edges in $H$ is less than $e(H)d^2 \le nd^3/2$,
the probability that any three such edges fall in the same $K_t$ copy of ${\cal D}$ is at most
\begin{equation}\label{e:4v}
d^3n \cdot \frac{t^2}{n^2} \le \frac{d^3t^2}{n}\;.
\end{equation}
Summing (\ref{e:bad},\ref{e:6v},\ref{e:5v},\ref{e:4v}), we obtain that the probability that $H_\sigma$ is atypical is at most
$$
\frac{8t^3d^2}{n}+ \frac{d^3t^4}{n} + \frac{d^3t^3}{n} + \frac{d^3t^2}{n} < \frac{11d^3t^4}{n}\;.
$$
Hence the number of atypical labeled copies of $H$ in $K_n$ is at most $11d^3t^4(n-1)!$.
\qed

\begin{definition}\label{d:pair-types}(consistency, inconsistency, disjointness)
A pair of edges of an orientation is called {\em consistent} if they share a common vertex and that vertex is the tail of one of them and the head of the other.
A pair of edges of an orientation is {\em inconsistent} if they share a common vertex and that vertex is the tail of both of them or the head of both of them.
A pair of edges of an orientation that is neither consistent nor inconsistent is {\em disjoint} (they do not share a common vertex).
\end{definition}
Observe that $plus(H)$ is the number of consistent pairs while $minus(H)$ is the number of inconsistent pairs.

Let $T$ be a tournament on vertex set $[t]$. For a permutation $\sigma \in S_t$ we let $T_\sigma$ denote the automorphism of $T$ obtained by the
permutation $\sigma$ on the vertex set $[t]$.
\begin{lemma}\label{l:prob-kt}
Let $t \ge 3$ be odd. Let $R$ be a regular tournament on vertex set $[t]$.
Consider the symmetric probability space whose elements are all the $R_\sigma$ for $\sigma \in S_t$.
Then, for any ordered triple $(x,y,z)$ of distinct vertices from $[t]$, the probability that the edge with endpoints $\{x,y\}$ is consistent with the edge with endpoints
$\{y,z\}$ is $(t-1)/(2t-4)$. Equivalently, the probability that this pair of edges is inconsistent is $(t-3)/(2t-4)$.
The probability that the triple $x,y,z$ induces a $C_3$ is $\frac{t+1}{4(t-2)}$. Equivalently, the probability that it induces a $T_3$ is $\frac{3(t-3)}{4(t-2)}$.
\end{lemma}
\Proof
Suppose, without loss of generality that $(x,y)$ is an edge of $R_\sigma$. Since $R$ is a regular tournament, $y$ is a head of $(t-1)/2$ other edges connecting $y$ to the remaining $t-2$ vertices.
Thus, the probability that $(y,z) \in R_\sigma$ is $((t-1)/2)/(t-2)=(t-1)/(2t-4)$. The complementary event of being inconsistent has probability $1-(t-1)/(2t-4)=(t-3)/(2t-4)$.
A regular tournament with $t$ vertices has exactly $t(t-1)(t-3)/8$ copies of $T_3$. Hence the probability that $x,y,z$ induce a $T_3$ is
$\frac{3(t-3)}{4(t-2)}$. The complementary event of being a $C_3$ has probability $\frac{t+1}{4(t-2)}$.
\qed

Finally, we need the following lemma.
\begin{lemma}\label{l:choice-t}
For all $\epsilon > 0$, if  $t$ is a sufficiently large integer, then
$$
\frac{t+1}{t-2} \ge \left(\frac{t-1}{t-2}\right)^{3-\epsilon}\;.
$$
For all $K > 0$, if $\rho > 0$ is sufficiently small, then
$$
\left(1-\rho^2\right)^K(1+\rho) \ge 1+\frac{\rho}{2}\;.
$$
\end{lemma}
\Proof
The first inequality is equivalent to 
$$
(t+1)(t-2)^2 \ge (t-1)^3 \left(\frac{t-2}{t-1}\right)^{\epsilon}
$$
which, in turn, is equivalent to
$$
1+\frac{1}{t-2} \ge \left(1 + \frac{3t-5}{t^3-3t^2+4}\right)^{\frac{1}{\epsilon}}
$$
which clearly holds for all sufficiently large $t$ as a function of $\epsilon$.
For the second claim, we may clearly assume that $K$ is an integer, and the left hand side is
$$
\left(1-\rho^2\right)^K(1+\rho)  =  1+\rho-K\rho^2 -O(\rho^3)
$$
which is larger than $1+\rho/2$ for all sufficiently small $\rho$ as a function of $K$.
\qed

\section{Proof of Theorem \ref{t:main}}

A pair of edges with a common endpoint is {\em induced} if the two non-common endpoints are not connected.
For an orientation $H$, let $c(H)$ denote the number of induced consistent pairs of edges of $H$ and denote by $i(H)$ the number of induced inconsistent pairs of edges of $H$.
Let $f(H)$ denote the number of directed triangles (copies of $C_3$) in $H$ and let $g(H)$ denote the number of acyclic triangles (copies of $T_3$) in $H$.

It will be more convenient to prove the following theorem which implies Theorem \ref{t:main}.
\begin{theorem}\label{t:special}
For any $\epsilon > 0$ and any positive integer $k$, there exist $\delta=\delta(\epsilon,k)$ and $N=N(\epsilon,k)$ such that the following holds.
Suppose $n > N$ is odd and $H$ is an orientation with $n$ vertices, $\Delta(H) \le k$ and $c(H)+(3-\epsilon)f(H)-i(H)-g(H) \ge \epsilon n$.
Then $R(H) \ge (1+\delta)n!/2^{e(H)}$.
\end{theorem}

\begin{lemma}\label{l:implies}
Theorem \ref{t:special} implies Theorem \ref{t:main}.
\end{lemma}
\Proof
Let $H$ be an orientation. Note that each $C_3$ in $H$ consists of three directed paths of length $2$, while each $T_3$ in $H$ consists of a single directed path of length $2$
and two anti-directed paths of length $2$. Thus,
\begin{equation}\label{e:eq}
plus(H) = 3f(H)+c(H)+g(H)~~,~~minus(H)=2g(H)+i(H)\;.
\end{equation}
In particular $plus(H)-minus(H)=c(H)+3f(H)-i(H)-g(H)$.
Let $\epsilon > 0$ and $k$ be a positive integer.
Let $\gamma =\epsilon/(k^2+1)$ and suppose that $n > N(\gamma,k)$ is odd where $N$ in the constant from Theorem \ref{t:special}.
Now, suppose that $H$ is $(\epsilon,k)$-consistent with $n$ vertices. Thus, in particular, $plus(H)-minus(H) \ge \epsilon n$.
Observe that since $\delta(H) \le k$, then the number of $C_3$ in $G$ is smaller than $k^2 n$, that is, $f(H) \le k^2 n$.
Now,
$$
c(H)+(3-\gamma)f(H)-i(H)-g(H) = plus(H)-minus(H)-\gamma f(H) \ge \epsilon n -\gamma k^2 n \ge \gamma n.\;
$$
Hence, by Theorem \ref{t:special} with $\delta=\delta(\gamma,k)$, 
there are regular tournaments with $n$ vertices that contain at least $(1+\delta)n!/2^{e(H)}$ labeled copies of $H$.
\qed

{\bf Proof of Theorem \ref{t:special}.}
Let $\rho=\frac{1}{t-2}$ where $t$ is the least positive integer satisfying the following three inequalities:
\begin{eqnarray*}
\frac{t+1}{t-2} & \ge & \left(\frac{t-1}{t-2}\right)^{3-\epsilon},\\
\left(1-\rho^2\right)^{\frac{2k^2}{\epsilon}}(1+\rho) & \ge & 1+\frac{\rho}{2},\\
t & \ge &\frac{2}{\epsilon}\;.
\end{eqnarray*}
By Lemma \ref{l:choice-t}, $t$ is well-defined. Define $\delta=\rho/4=\frac{1}{4t-8}$.
$N$ will be chosen sufficiently large as a function of $\epsilon,k,t$ (and thus, as a function of $\epsilon,k$) so that all claimed inequalities involving
all $n > N$ hold throughout the proof.

Suppose that $n > N$ is odd and let  ${\cal D}$ be an adjusted $t$-decomposition of $K_n$ which exists by Lemma \ref{l:dec} and our choice of constants.
We consider the following probability space of regular tournaments on vertex set $[n]$.
Let $R$ be any fixed regular tournament on the set of vertices $[t]$. Let $R^*$ be any fixed regular tournament
on the set of vertices $[2t-1]$. For each element $Z \in {\cal D}$ perform the following. If $Z$ is isomorphic to $K_t$, then orient it by placing a random copy of $R$ on it.
More formally, suppose $Z$ consists of vertices $\{v_1,\ldots,v_t\}$. Take a random permutation $\sigma \in S_t$ and for each edge $v_i v_j$
orient it in the same direction of the edge of $R$ connecting $\sigma(i)$ and $\sigma(j)$. If $Z$ is isomorphic to $K_{2t-1}$ (recall, there are only at most $t-1$ such $Z$ in ${\cal D}$),
then place a random copy of $R^*$ on it. If $Z$ is a $4$-cycle, then orient the cycle as a directed cycle
in one of the two possible directions, each direction equally likely. If $Z$ is a $3$-cycle, then orient the cycle as a directed cycle
in one of the two possible directions, each direction equally likely. This regularly orients each element of ${\cal D}$ (the $|{\cal D}|$ random orientations are independent of each other)
and so the obtained random tournament is regular on vertex set $[n]$. Denote the obtained random regular tournament by $T$.

Let $\hat{H}$ denote the underlying undirected graph of $H$.
Recall that there are $n!$ labeled copies of $\hat{H}$ in $\hat{T}=K_n$ and we denote them by $\hat{H}_\pi$ for $\pi \in S_n$.
Each $\hat{H}_\pi$ corresponds to a labeling $H_\pi$ of $H$. However, while $\hat{H}_\pi$ is trivially a labeled copy of $\hat{H}$ in $K_n$,
$H_\pi$ is not necessarily a labeled copy of $H$ in $T$, since the direction of the edges may not agree with the direction of the edges in $T$.
Let ${\cal L}$ denote the typical labeled copies of $\hat{H}$ in $K_n$ with respect to ${\cal D}$ (recall Definition \ref{d:atypical}).
By Lemma \ref{l:typical} we know that $|{\cal L}| \ge n! - 11k^3t^4(n-1)!$. We will prove that the expected number of typical labeled copies of $H$ in $T$
is already $(1+\delta)n!/2^{e(H)}$. In particular, this implies that there exists a regular $n$-vertex tournament for which the number of labeled copies of $H$ in it is at least
$(1+\delta)n!/2^{e(H)}$.

For a labeled copy $\hat{H}_\pi$, let $C_\pi$ denote the number of induced consistent pairs of edges of $H_\pi$ that fall into the same element of ${\cal D}$,
let $I_\pi$ denote the number of induced inconsistent pairs of edges of $H_\pi$ that fall into the same element of ${\cal D}$, let $F_\pi$ denote the number of $C_3$ of $H_\pi$ 
that fall into the same element of ${\cal D}$, and let $G_\pi$ denote denote the number of $T_3$ of $H_\pi$ that fall into the same element of ${\cal D}$.
Trivially, $C_\pi \le c(H)$, $I_\pi \le i(H)$, $F_\pi \le f(H)$ and $G_\pi \le g(H)$ but they will usually be much smaller.

Now, if $\hat{H}_\pi$ is typical (namely, $\hat{H}_\pi \in {\cal L}$), observe that if two induced consistent or inconsistent pairs of edges fall into the same element of ${\cal D}$, then this
element must be a $K_t$. Similarly, if the edges of a copy of $C_3$ or the edges of a copy of $T_3$ (namely, the edges of a triangle of $\hat{H}_\pi$) fall into the same element of ${\cal D}$,
then this element must be a $K_t$.
Let $X_\pi$ be the indicator random variable which equals $1$ if $H_\pi$ is a labeled copy of $H$ in $T$ (recall: the random object here is $T$).
Notice first an easy case: if no two edges of $H_\pi$ are in the same element of ${\cal D}$, then trivially $\Pr[X_\pi = 1] = 2^{-e(H)}$, since recall that each element of ${\cal D}$ is oriented
independently, and any single edge has equal probability $\frac{1}{2}$ in both directions. The following claim determines $\Pr[X_\pi = 1]$ for any typical copy $\hat{H}_\pi \in {\cal L}$
in terms of $C_\pi, I_\pi, F_\pi, G_\pi$.

\begin{claim}\label{cl:1}
Suppose $\hat{H}_\pi \in {\cal L}$, then
$$
\Pr[X_\pi = 1] = \left( \frac{t-1}{4(t-2)}\right)^{C_\pi} \left( \frac{t-3}{4(t-2)}\right)^{I_\pi}\left( \frac{t+1}{8(t-2)}\right)^{F_\pi} \left( \frac{t-3}{8(t-2)}\right)^{G_\pi} 
\left( \frac{1}{2}\right)^{e(H)-2C_\pi - 2I_\pi-3F_\pi-3G_\pi}\;.
$$
\end{claim}
\Proof
We say that an element of ${\cal D}$ is {\em good} (with respect to $H_\pi$ and $T$) if it does not contain edges of $H_\pi$, or if all the edges it contains from $H_\pi$ agree with $T$ (have the same orientation as in $T$).
If an element of ${\cal D}$ contains a single edge of $H_\pi$, then it is good with probability $\frac{1}{2}$ since recall that
in $T$ each edge direction is equally likely. Since $\hat{H}_\pi$ is typical, no element of ${\cal D}$ contains three edges of $H_\pi$
unless these three edges are a triangle of $\hat{H}_\pi$.
So, we now consider an element of ${\cal D}$ that contains precisely two edges of $\hat{H}_\pi$, or a triangle of $\hat{H}_\pi$.
Again, by typicality, this element is a $K_t$.

Assume first that it contains precisely two edges and these two edges do not share a common endpoint.
Since $K_t$ was oriented in $T$ by taking a random permutation of $R$ (which is a regular tournament with $t$ vertices), the probability that both edges agree
with $T$ is $(1/2) \cdot (1/2)=1/4$, so the probability that this element of ${\cal D}$ is good is $1/4$.
Assume next that these two edges share a common endpoint. There are two cases here. Either these two edges in $H_\pi$ induce a consistent pair in $H_\pi$, or they induce an inconsistent
pair in $H_\pi$.
If they are consistent, say they are $(x,y)$ and $(y,z)$, then the probability that $(x,y)$ agrees with $T$ is $\frac{1}{2}$, and given that it agrees with $T$,
the probability that also $(y,z)$ agrees with $T$ is, by Lemma \ref{l:prob-kt}, $(t-1)/(2t-4)$. Thus, the probability that this element  of ${\cal D}$ is good is 
$\frac{t-1}{4(t-2)}$. If they are inconsistent, say they are $(x,y)$ and $(z,y)$, then the probability that $(x,y)$ agrees with $T$ is $\frac{1}{2}$, and given that it agrees with $T$,
the probability that also $(z,y)$ agrees with $T$ is, by Lemma \ref{l:prob-kt}, $(t-3)/(2t-4)$. Thus, the probability that this element  of ${\cal D}$ is good is 
$\frac{t-3}{4(t-2)}$.

Assume next that it contains three edges which form a triangle in $\hat{H}_\pi$. Then it is either a $C_3$ in $H_\pi$ or a $T_3$ in $H_\pi$.
Since $K_t$ was oriented in $T$ by taking a random permutation of $R$, the probability that the edges of a $C_3$ of $H_\pi$ falling in the same $K_t$ copy all agree
with $T$ is, by Lemma \ref{l:prob-kt}, precisely $(t+1)/(8(t-2))$ since by Lemma \ref{l:prob-kt}, the three vertices are a directed triangle with probability 
$(t+1)/(4(t-2))$, and the two possible directions of the triangle are equally likely.
Similarly, the probability that the edges of a $T_3$ of $H_\pi$ falling in the same $K_t$ copy all agree
with $T$ is, by Lemma \ref{l:prob-kt}, precisely $(t-3)/(8(t-2))$ since by Lemma \ref{l:prob-kt}, the three vertices are a $T_3$ with probability 
$3(t-3)/(4(t-2))$, and the six possible orientations of that $T_3$ are equally likely. Hence, the probability that the element of ${\cal D}$ is good in
the $C_3$ case is $(t+1)/(8(t-2))$ and is $\frac{t-3}{8(t-2)}$ in the $T_3$ case.

Since the goodness of each element of ${\cal D}$ is independent of the goodness of any other element, and since there are $C_\pi$ elements that contain induced consistent pairs,
$I_\pi$ elements that contain induced inconsistent pairs, $F_\pi$ elements that contain $C_3$'s and $G_\pi$ elements that contains $T_3$'s,
we obtain the expression in the statement of the claim for the probability that all elements of ${\cal D}$ are good.
\qed

It follows from Claim \ref{cl:1} that the expected number of typical labeled copies of $H$ in $T$ is at least
\begin{equation}\label{e:exp}
\sum_{\hat{H}_\pi \in {\cal L}} \left( \frac{t-1}{4(t-2)}\right)^{C_\pi} \left( \frac{t-3}{4(t-2)}\right)^{I_\pi}
\left( \frac{t+1}{8(t-2)}\right)^{F_\pi} \left( \frac{t-3}{8(t-2)}\right)^{G_\pi}\left( \frac{1}{2}\right)^{e(H)-2C_\pi - 2I_\pi-3F_\pi-3G_\pi}\;.
\end{equation}
It hence remains to show that (\ref{e:exp}) is at least $(1+\delta)n!/2^{e(H)}$.

Let $C_{avg}$ denote the average of $C_\pi$ ranging over all $\hat{H}_\pi \in {\cal L}$ and let $I_{avg}$ denote the average of $I_\pi$ ranging over all $\hat{H}_\pi \in {\cal L}$.
Similarly define $F_{avg}$ and $G_{avg}$.
Let $C^*_{avg}$ denote the average of $C_\pi$ ranging over all $n!$ labeled copies $\hat{H}_\pi$ and 
let $I^*_{avg}$ denote the average of $I_\pi$ ranging over all $n!$ labeled copies $\hat{H}_\pi$
(so the difference between $C^*_{avg}$ and $C_{avg}$ is that the former averages over all labeled copies while the latter averages only over typical labeled copies
and similarly for $I^*_{avg}$ and $I_{avg}$). Similarly define $F^*_{avg}$ and $G^*_{avg}$.

Let us first use $C_{avg},I_{avg},F_{avg},G_{avg}$ to obtain a lower bound for (\ref{e:exp}).
\begin{claim}\label{cl:2}
\begin{eqnarray}
& \sum_{\hat{H}_\pi \in {\cal L}} \left( \frac{t-1}{4(t-2)}\right)^{C_\pi} \left( \frac{t-3}{4(t-2)}\right)^{I_\pi}
\left( \frac{t+1}{8(t-2)}\right)^{F_\pi} \left( \frac{t-3}{8(t-2)}\right)^{G_\pi}\left( \frac{1}{2}\right)^{e(H)-2C_\pi - 2I_\pi-3F_\pi-3G_\pi} \nonumber\\
\ge & |{\cal L}| \cdot \left( \frac{t-1}{4(t-2)}\right)^{C_{avg}} \left( \frac{t-3}{4(t-2)}\right)^{I_{avg}}
\left( \frac{t+1}{8(t-2)}\right)^{F_{avg}} \left( \frac{t-3}{8(t-2)}\right)^{G_{avg}}
 \left( \frac{1}{2}\right)^{e(H)-2C_{avg} - 2I_{avg}-3F_{avg} - 3G_{avg}}\;.\label{e:avg}
\end{eqnarray}
\end{claim}
\Proof
The inequality follows immediately from the inequality of arithmetic and geometric means.
\qed

So, our remaining task is to lower bound (\ref{e:avg}), i.e. the right hand side of the inequality in Claim \ref{cl:2}.
We will do so gradually by first proving good estimates for $C^*_{avg},I^*_{avg},F^*_{avg},G^*_{avg}$ and then proving that 
$C_{avg}$ is very close to $C^*_{avg}$ and $I_{avg}$ is very close to $I^*_{avg}$ and respectively for $F_{avg}$ and $G_{avg}$.

Let us consider first a case where $C^*_{avg},I^*_{avg},F^*_{avg},G^*_{avg}$ are straightforward to compute.
This is the case where all elements of ${\cal D}$ are isomorphic to $K_t$ (namely when the adjusted $t$-decomposition is actually a $K_t$-decomposition).
Instead of counting the labeled copies is which two edges (or three edges of a triangle) appear in the same element of ${\cal D}$, it will be more convenient (though equivalent) to consider
the labeled copies as a symmetric probability space of size $n!$ where each labeled copy is equally likely.
Consider a pair of edges of $H$ with a common endpoint (consistent, inconsistent, or inducing a triangle). Say the endpoints of the first edge are $u,v$ and the endpoints of the second edge is
$v,w$. The probability that a randomly chosen labeled copy of $\hat{H}$ will place $u,v,w$ all in the same copy of ${\cal D}$ is precisely $(t-2)/(n-2)$.
Since there are $c(H)$ consistent pairs in $H$, we have by linearity of expectation that $C^*_{avg} = c(H)(t-2)/(n-2)$ and similarly $I^*_{avg} = i(H)(t-2)/(n-2)$,
$F^*_{avg} = f(H)(t-2)/(n-2)$,$G^*_{avg} = g(H)(t-2)/(n-2)$.
The following claim estimates $C^*_{avg},I^*_{avg},F^*_{avg},G^*_{avg}$ in the case that ${\cal D}$ is an arbitrary adjusted $t$-decomposition.
\begin{claim}\label{cl:avg*}
For a pair of edges of $\hat{H}$ with a common endpoint, let $p$ denote the probability that a labeled copy of $\hat{H}$ chosen at random will place both edges in the same element of ${\cal D}$.
Then,
$$
\left|p-\frac{t-2}{n-2}\right| \le \frac{3t^2}{n^2}\;.
$$
\end{claim}
\Proof
Let ${\cal D}'$ be the subset of elements of ${\cal D}$ that are not isomorphic to $K_t$. By the definition of adjusted $t$-decomposition, all elements of ${\cal D}'$ together
contain at most $n(t-3)/2+4(t-3)+(t-1)^2(2t-1) < nt$ edges.
 
Consider a pair of edges $e,f$ of $\hat{H}$ with a common endpoint. Say the endpoints of $e$ are $u,v$ and the endpoints of $f$ are $v,w$.
Let $p$ denote the probability that a random labeled copy of $\hat{H}$ will place both $e$ and $f$ in the same element of ${\cal D}$.
Let $p_1$ denote the probability that $e$ is placed in an element of ${\cal D'}$.
Let $q_1$ denote the probability that $f$ is placed in the same element as $e$ given that $e$ is placed in an element of ${\cal D'}$.
Let $q_2$ denote the probability that $f$ is placed in the same element as $e$ given that $e$ is placed in a $K_t$ element of ${\cal D}$.
Thus, $p=p_1q_1+(1-p_1)q_2$.

Since all elements of ${\cal D}'$ together contain at most $nt$ edges, we have that $p_1 \le nt/\binom{n}{2}=2t/(n-1)$.
Since each element of ${\cal D}'$ is either a triangle or a cycle of length $4$ or a $K_{2t-1}$, the probability that vertex $w$ (and hence $f$) is also placed in the same element of ${\cal D}'$ in which $e$ is placed is
at most $(2t-3)/(n-2)$. Hence, $q_1 \le (2t-3)/(n-2)$. Given that $e$ is placed in a $K_t$-element, the probability that vertex $w$ (and hence $f$) is also placed in that element is
$q_2=(t-2)/(n-2)$.
Hence,
$$
\left|p-\frac{t-2}{n-2}\right| = \left|p-q_2\right| = p_1|q_1-q_2| \le \frac{2t}{n-1} \cdot \frac{t-1}{n-2} \le \frac{3t^2}{n^2}\;.
$$
\qed

Now, since $C^*_{avg} = c(H)p$, $I^*_{avg} = i(H)p$, $F^*_{avg} = f(H)p$, $G^*_{avg} = g(H)p$ where $p$, as in the last claim, is the probability that two edges sharing an endpoint
are placed in the same element of ${\cal D}$, we immediately obtain the following.
\begin{corollary}\label{c:avg-est}
$$
\left|C^*_{avg} - c(H)\frac{t-2}{n-2}\right| \le c(H) \frac{3t^2}{n^2}~~,~~\left|I^*_{avg} - i(H)\frac{t-2}{n-2}\right| \le i(H) \frac{3t^2}{n^2}\;.
$$
$$
\left|F^*_{avg} - f(H)\frac{t-2}{n-2}\right| \le f(H) \frac{3t^2}{n^2}~~,~~\left|G^*_{avg} - g(H)\frac{t-2}{n-2}\right| \le g(H) \frac{3t^2}{n^2}\;.
$$

\qed
\end{corollary}

\noindent
We next show that $C^*_{avg}$ and $C_{avg}$ are not far apart and similarly for $(I^*_{avg},I_{avg})$, $(F^*_{avg},F_{avg})$, $(G^*_{avg},G_{avg})$.
\begin{claim}\label{cl:avg}
$$
c(H)\left(\frac{t-2}{n-2}-\frac{3t^2}{n^2}\right)-\frac{133k^3t^4d^2}{\sqrt{n}} \le C_{avg} \le c(H) \left(\frac{t-2}{n-2}+\frac{23k^3t^5}{n^2}\right)\;.
$$
$$
i(H)\left(\frac{t-2}{n-2}-\frac{3t^2}{n^2}\right)-\frac{133k^3t^4d^2}{\sqrt{n}} \le I_{avg} \le i(H) \left(\frac{t-2}{n-2}+\frac{23k^3t^5}{n^2}\right)\;.
$$
$$
f(H)\left(\frac{t-2}{n-2}-\frac{3t^2}{n^2}\right)-\frac{133k^3t^4d^2}{\sqrt{n}} \le F_{avg} \le f(H) \left(\frac{t-2}{n-2}+\frac{23k^3t^5}{n^2}\right)\;.
$$
$$
g(H)\left(\frac{t-2}{n-2}-\frac{3t^2}{n^2}\right)-\frac{133k^3t^4d^2}{\sqrt{n}} \le G_{avg} \le g(H) \left(\frac{t-2}{n-2}+\frac{23k^3t^5}{n^2}\right)\;.
$$
\end{claim}
\Proof
We prove for $C_{avg}$. The proof for the others is identical.
Let $C(H)$ be the set of triplets of vertices of $H$ that induce a consistent pair in $H$, thus $|C(H)|=c(H)$. A subset $S \subset C(H)$ is {\em independent} if any two elements of $S$
are disjoint (so their union is $6$ vertices).
For an independent subset $S \subset C(H)$ we say that a labeled copy of $\hat{H}$ is {\em $S$-bad} if for each triplet in $S$, all three vertices of the triplet appear in the same
element of ${\cal D}$. Our goal is to prove that if $|S|$ is large, the probability of being $S$-bad is very small.

Fix some $S \subset C(H)$ with $s=|S|$ and suppose that $S=\{\{x_i,y_i,z_i\} ~|~ i=1,\ldots,s\}$.
What is the probability that a random labeled copy of $\hat{H}$ is $S$-bad? Recall that each element of ${\cal D}$ is either a $K_t$, a $C_3$, a $C_4$ or a $K_{2t-1}$.
Consider first the triplet $\{x_1,y_1,z_1\}$. The probability that all three fall in the same element of ${\cal D}$ is at most $(2t-3)/(n-2)$ (it is actually smaller since most elements
of ${\cal D}$ are $K_t$ so the probability is very close to $(t-2)/(n-2)$ but we don't want to bother with a case analysis).
Given that $\{x_i,y_i,z_i\}$ each fall in the same element of ${\cal D}$ for $i=1,\ldots,r-1$ (distinct triplets may or may not fall in the same element of ${\cal D}$) what is the probability
that also
$\{x_r,y_r,z_r\}$ fall in the same element of ${\cal D}$? Suppose we are given the information to which element of ${\cal D}$ each of the $3r-1$ vertices of
$x_1,y_1,z_1,x_2,y_2,z_2\ldots,x_{r-1},y_{r-1},z_{r-1},x_r,y_r$ 
fall. The probability that $z_r$ also falls in the copy to which $x_r$ and $y_r$ fell is thus at most $(2t-3)/(n-3r+1)$.
Hence, the probability that a random labeled copy of $\hat{H}$ is $S$-bad is at most
$$
\left(\frac{2t-3}{n-3s+1}\right)^s\;.
$$
Recall that $\Delta(H) \le k$ and therefore, trivially $c(H) < k^2 n$ (and the same trivial upper bound holds for $i(H),f(H),g(H)$). Let us use, say, $s=\lfloor \sqrt{n} \rfloor$ for the remainder of the proof
and say that a labeled copy of $\hat{H}$ is {\em bad} if it is $S$-bad for
some independent $S \subset C(H)$ with $|S|=s=\lfloor \sqrt{n} \rfloor$. Thus, the probability that a random labeled copy of $\hat{H}$ is bad is at most
\begin{equation}\label{e:vbad}
\binom{c(H)}{s} \left(\frac{2t-3}{n-3s+1}\right)^s \le \left(\frac{e k^2n}{s}\right)^s \left(\frac{2t-3}{n-3s+1}\right)^s \ll \frac{1}{n^2}
\end{equation}
where we have only used here that $n$ is sufficiently large as a function of $t$ and $k$.
In other words, we have proved hat the number of bad labeled copies of $\hat{H}$ is at most $n!/n^2 \le (n-2)!$.

We will say that a labeled copy of $\hat{H}$ is {\em very bad} if at least $6k^2s$ elements of $C(H)$ each appear together in the same element of ${\cal D}$.
So the difference between being bad and very bad is that in the latter we don't insist that the triplets are independent.
We claim that if a labeled copy is very bad, then it is also bad. Indeed, consider the graph $W$ whose vertices are the elements of $C(H)$, and two vertices of $W$ are connected if they are not
disjoint (i.e. they span together less than $6$ vertices). Since $\Delta(H) \le k$, each vertex of $H$ appears in less than
$2k^2$ triplets. Hence, the maximum degree of $W$ is less than $6k^2$. It follows that any set of $6k^2s$ elements of $C(H)$ has an independent set of size
at least $s$, hence if a labeled copy is very bad, then it is also bad, and thus, by (\ref{e:vbad}), there are at most $(n-2)!$ very bad labeled copies of $\hat{H}$.

We are now ready to compute our upper and lower bounds for $C_{avg}$.
Recall that by their definitions,
$$
\sum_{\pi \in S_n} C_\pi  = n! C^*_{avg}~~,~~ \sum_{\hat{H}_\pi \in{\cal L}} C_\pi  = |{\cal L}| C_{avg}\;.
$$
The upper bound follows rather simply from Lemma \ref{l:typical} which asserts that $|{\cal L}| \ge n! - 11k^3t^4(n-1)!$.
Indeed,
$$
C_{avg}|{\cal L}| = n! C^*_{avg} - \sum_{\hat{H}_\pi \notin{\cal L}} C_\pi \le n! C^*_{avg}\;.
$$
So,
\begin{eqnarray*}
C_{avg} & \le & \frac{n! C^*_{avg}}{|L|}  \\
& = & C^*_{avg} + \frac{(n!-|{\cal L}|)C^*_{avg}}{|{\cal L}|} \\
& \le & C^*_{avg} + \frac{11k^3t^4(n-1)! C^*_{avg}}{|{\cal L}|} \\
& \le & C^*_{avg} + \frac{22k^3t^4 C^*_{avg}}{n} \\
& = & C^*_{avg}\left(1+\frac{22k^3t^4 }{n}\right)\\
& \le & c(H)\left(\frac{t-2}{n-2}+\frac{3t^2}{n^2}\right)\left(1+\frac{22k^3t^4 }{n}\right)\\
& \le & c(H) \left(\frac{t-2}{n-2}+\frac{23k^3t^5}{n^2}\right)\;.
\end{eqnarray*}
Notice that we have used Claim \ref{c:avg-est} which asserts that $C^*_{avg} \le c(H)\left(\frac{t-2}{n-2}+\frac{3t^2}{n^2}\right)$.

For the lower bound, we need first to bound $\sum_{\hat{H}_\pi \notin{\cal L}} C_\pi$.
We know that the number of $\hat{H}_\pi \notin{\cal L}$ is at most $11k^3t^4(n-1)!$. By the first part of the proof, at most $(n-2)!$ of them are very bad,
and for them we will use the trivial bound $C_\pi \le c(H) \le k^2n$. The others are not very bad, so for them we have $C_\pi \le 6k^2\sqrt{n}$.
Thus,
$$
\sum_{\hat{H}_\pi \notin{\cal L}} C_\pi \le 11k^3t^4(n-1)! \cdot 6k^2\sqrt{n} + (n-2)!k^2 n\;.
$$
Now,
\begin{eqnarray*}
C_{avg} & = & \frac{n! C^*_{avg} - \sum_{\hat{H}_\pi \notin{\cal L}} C_\pi}{|{\cal L}|}  \\
& \ge & C^*_{avg} - \frac{\sum_{\hat{H}_\pi \notin{\cal L}} C_\pi}{|{\cal L}|} \\
& \ge & C^*_{avg} - \frac{2\sum_{\hat{H}_\pi \notin{\cal L}} C_\pi}{n!} \\
& \ge & C^*_{avg} - \frac{2 \cdot 11k^3t^4\cdot 6k^2}{\sqrt{n}} - \frac{2k^2}{n-1}\\
& \ge & C^*_{avg} - \frac{133k^5t^4}{\sqrt{n}} \\
& \ge & c(H)\left(\frac{t-2}{n-2}-\frac{3t^2}{n^2}\right)-\frac{133k^5t^4}{\sqrt{n}}\;.
\end{eqnarray*}
Notice that we have used Claim \ref{c:avg-est} which asserts that $C^*_{avg} \ge c(H)\left(\frac{t-2}{n-2}-\frac{3t^2}{n^2}\right)$.
\qed

\begin{corollary}\label{c:diff}
$C_{avg} + (3-\epsilon)F_{avg} - I_{avg} - G_{avg} \ge \frac{\epsilon}{2}t$.
\end{corollary}
\Proof
By Claim \ref{cl:avg},
$$
C_{avg} + (3-\epsilon)F_{avg} - I_{avg} - G_{avg} \ge
$$
$$
(c(H)+(3-\epsilon)f(H))\left(\frac{t-2}{n-2}-\frac{3t^2}{n^2}\right)-4\frac{133k^5t^4}{\sqrt{n}} - (i(H)+g(H)) \left(\frac{t-2}{n-2}+\frac{23k^3t^5}{n^2}\right)\;.
$$
By the assumption of Theorem \ref{t:special}, $c(H)+(3-\epsilon)f(H)-i(H)-g(H) \ge \epsilon n$. Also, trivially, each of $c(H),i(H),f(H),g(H)$ is at most  $k^2n$
so the last inequality is at least
$$
\epsilon n \frac{t-2}{n-2}-O\left(\frac{1}{\sqrt{n}}\right) \ge \frac{\epsilon}{2}t
$$
for all $n$ sufficiently large as a function of $\epsilon,k,t$.
\qed

We are now ready to prove our final claim which completes the proof of Theorem \ref{t:special}.
\begin{claim}\label{cl:special}
$$
|{\cal L}| \cdot \left( \frac{t-1}{4(t-2)}\right)^{C_{avg}} \left( \frac{t-3}{4(t-2)}\right)^{I_{avg}}
\left( \frac{t+1}{8(t-2)}\right)^{F_{avg}} \left( \frac{t-3}{8(t-2)}\right)^{G_{avg}}
 \left( \frac{1}{2}\right)^{e(H)-2C_{avg} - 2I_{avg}-3F_{avg} - 3G_{avg}} 
$$
$$
 \ge n!(1+\delta)\left( \frac{1}{2}\right)^{e(H)}\;.
$$
\end{claim}
\Proof
First observe that since $|{\cal L}| \ge n!-11k^3t^4(n-1)!$, then for all $n$ sufficiently large as a function of $\delta,t,k$ we have
$|{\cal L}| \ge (1+\delta)n!/(1+2\delta)$. Hence it suffices to prove that
$$
\left( \frac{t-1}{4(t-2)}\right)^{C_{avg}} \left( \frac{t-3}{4(t-2)}\right)^{I_{avg}}
\left( \frac{t+1}{8(t-2)}\right)^{F_{avg}} \left( \frac{t-3}{8(t-2)}\right)^{G_{avg}}
 \left( \frac{1}{2}\right)^{e(H)-2C_{avg} - 2I_{avg}-3F_{avg} - 3G_{avg}} 
$$
$$
 \ge (1+2\delta)\left( \frac{1}{2}\right)^{e(H)}\;.
$$
which is equivalent to proving that
$$
\left( \frac{t-1}{t-2}\right)^{C_{avg}} \left( \frac{t-3}{t-2}\right)^{I_{avg}}
\left( \frac{t+1}{t-2}\right)^{F_{avg}} \left( \frac{t-3}{t-2}\right)^{G_{avg}}
 \ge 1+2\delta\;.
$$
Recall that by our choice of $t$,
$$
\frac{t+1}{t-2} \ge \left(\frac{t-1}{t-2}\right)^{3-\epsilon}\;.
$$
Hence it suffices to prove that
$$
\left( \frac{t-1}{t-2}\right)^{C_{avg}+(3-\epsilon)F_{avg}} \left( \frac{t-3}{t-2}\right)^{I_{avg}+G_{avg}} \ge 1+2\delta\;.
$$
Now, by Claim \ref{cl:avg} and the trivial bound $i(H)+g(H) \le k^2 n$ we have
\begin{equation}\label{e:k2t}
I_{avg}+G_{avg} \le (i(H)+g(H)) \left(\frac{t-2}{n-2}+\frac{23k^3t^5}{n^2}\right) \le (i(H)+g(H)) \frac{t-1}{n-2} \le k^2 t\;.
\end{equation}
Hence, by Corollary \ref{c:diff}
\begin{eqnarray*}
& & \left( \frac{t-1}{t-2}\right)^{C_{avg}+(3-\epsilon)F_{avg}} \left( \frac{t-3}{t-2}\right)^{I_{avg}+G_{avg}}\\
& = &\left( 1-\frac{1}{(t-2)^2}\right)^{I_{avg}+G_{avg}}\left(\frac{t-1}{t-2}\right)^{C_{avg}+(3-\epsilon)F_{avg}-I_{avg}-G_{avg}}\\
& \ge & \left( 1-\frac{1}{(t-2)^2}\right)^{I_{avg}+G_{avg}} \left(\frac{t-1}{t-2}\right)^{\frac{\epsilon}{2}t}\\
& \ge & \left( 1-\frac{1}{(t-2)^2}\right)^{k^2t } \left(\frac{t-1}{t-2}\right)^{\frac{\epsilon}{2}t}
\end{eqnarray*}
So it remains to prove that
$$
\left( 1-\frac{1}{(t-2)^2}\right)^{k^2t } \left(1+\frac{1}{t-2}\right)^{\frac{\epsilon}{2}t} \ge 1+2\delta
$$
which is equivalent to proving
$$
\left( 1-\frac{1}{(t-2)^2}\right)^{\frac{2k^2}{\epsilon}} \left(1+\frac{1}{t-2}\right) \ge \left(1+2\delta\right)^{\frac{2}{\epsilon t}}
$$
But recall that $\delta=\frac{1}{4(t-2)}$ and that $2/(\epsilon t) \le 1$ so the last inequality holds if the following inequality holds:
$$
\left( 1-\frac{1}{(t-2)^2}\right)^{\frac{2k^2}{\epsilon}} \left(1+\frac{1}{t-2}\right) \ge 1+\frac{1}{2t-4}
$$
and indeed the last inequality holds by our choice of $t$.
\qed

\vspace{10pt}
Claim \ref{cl:special} proves the claimed  lower bound for (\ref{e:avg}) and hence Theorem \ref{t:special} holds. By Lemma \ref{l:implies}, Theorem \ref{t:main} holds as well. \qed

\section{Balanced and Eulerian orientations}

{\bf Proof of Proposition \ref{p:euler}}.
Let $H$ be a balanced orientation with $n$ vertices, maximum degree $k$, and at least $(0.5+\epsilon)n$ edges.
We first notice that the underlying graph $\hat{H}$ has at least $2\epsilon n/k$ vertices with degree at least $2$.
Indeed, otherwise the sum of the degrees of the vertices of ${\hat H}$ would have been less than $n+2\epsilon n = (1+2\epsilon)n$, contradicting the fact that is has at least
$(0.5+\epsilon)n$ edges.

Since $H$ is balanced, we have that for any vertex $v$ whose degree in $H$ is at least $2$,
\begin{equation}\label{e:v-contrib}
d^+(v)d^{-}(v) - \left(\binom{d^+(v)}{2}+\binom{d^{-}(v)}{2}\right) \ge 1\;.
\end{equation}
It follows that
$$
plus(H)-minus(H) \ge 2\epsilon n/k\;
$$
which means that $H$ is a $(2\epsilon/k,k)$-consistent orientation.
If $H$ is an Eulerian orientation, then (\ref{e:v-contrib}) holds for each $v \in V(H)$ so $plus(H)-minus(H) \ge n$ which means that $H$ is a $(1,k)$-consistent orientation.
\qed

\vspace{10pt}

\noindent
{\bf Proof of Theorem \ref{t:k-reg}}.
We show that Theorem \ref{t:k-reg} follows from a careful optimization of the inequalities in the proof of Theorem \ref{t:special}.
We first observe that if $H$ is a $k$-regular orientation, then $plus(H)=k^2n$ while $minus(H)=k(k-1)n$ and recall that using the notation
in Section $3$ we have by (\ref{e:eq}) that $plus(H)=3f(H)+c(H)+g(H)=nk^2$ while $minus(H)=2g(H)+i(H)=k(k-1)n$.

Recalling Claim \ref{cl:2}, our goal is to prove that for every $\gamma > 0$, and for all sufficiently large $t$,
$$
|{\cal L}| \cdot \left( \frac{t-1}{4(t-2)}\right)^{C_{avg}} \left( \frac{t-3}{4(t-2)}\right)^{I_{avg}}
\left( \frac{t+1}{8(t-2)}\right)^{F_{avg}} \left( \frac{t-3}{8(t-2)}\right)^{G_{avg}}
 \left( \frac{1}{2}\right)^{e(H)-2C_{avg} - 2I_{avg}-3F_{avg} - 3G_{avg}}
$$
is at least  $(e^k-\gamma)n!/2^{nk}$ for all $n$ sufficiently large.
As $e(H)=nk$ this is equivalent to proving that for every $\gamma > 0$, and for all sufficiently large $t$,
$$
|{\cal L}| \left( \frac{t-1}{t-2}\right)^{C_{avg}} \left( \frac{t-3}{t-2}\right)^{I_{avg}}
\left( \frac{t+1}{t-2}\right)^{F_{avg}} \left( \frac{t-3}{t-2}\right)^{G_{avg}} \ge (e^k-\gamma)n!
$$
holds for all $n$ sufficiently large.

By Lemma \ref{l:typical}, $|{\cal L}| \ge n!-11k^3t^4(n-1)!$ so it is equivalent to prove that for all $t$ sufficiently large
$$
\left( \frac{t-1}{t-2}\right)^{C_{avg}} \left( \frac{t-3}{t-2}\right)^{I_{avg}}
\left( \frac{t+1}{t-2}\right)^{F_{avg}} \left( \frac{t-3}{t-2}\right)^{G_{avg}} \ge e^k-\frac{\gamma}{2}
$$
holds for all $n$ sufficiently large.

Simple algebraic manipulation of the left hand side shows that it is equal to
$$
\left( \frac{t-1}{t-2}\right)^{C_{avg}+3F_{avg}-I_{avg}-G_{avg}}\left( 1-\frac{1}{(t-2)^2}\right)^{I_{avg}+G_{avg}}\left(1-\frac{3t-5}{(t-1)^3}\right)^{F_{avg}}\;.
$$
But now, from Claim \ref{cl:avg} and the trivial bound $i(H)+g(H) \le k^2 n$ and $f(H) \le k^2n$ we have as in (\ref{e:k2t})
that
$$
I_{avg}+G_{avg} \le k^2 t~~,~~F_{avg} \le k^2t\;.
$$
Also, by Claim \ref{cl:avg},
$$
C_{avg} + 3F_{avg} - I_{avg} - G_{avg} \ge
(c(H)+3f(H))\left(\frac{t-2}{n-2}-\frac{3t^2}{n^2}\right)-4\frac{133k^5t^4}{\sqrt{n}} - (i(H)+g(H)) \left(\frac{t-2}{n-2}+\frac{23k^3t^5}{n^2}\right)\;.
$$
Now, since $c(H)+3f(H)-i(H)-g(H) = kn$ and since each of $c(H),i(H),f(H),g(H)$ is at most  $k^2n$ we have that
$$
C_{avg} + 3F_{avg} - I_{avg} - G_{avg} \ge kn \frac{t-2}{n-2}-O\left(\frac{1}{\sqrt{n}}\right) \ge k(t-3)
$$
for all $n$ sufficiently large as a function of $k,t$.
It thus remains to show that
$$
\left( \frac{t-1}{t-2}\right)^{k(t-3)}\left( 1-\frac{1}{(t-2)^2}\right)^{k^2t}\left(1-\frac{3t-5}{(t-1)^3}\right)^{k^2t} \ge e^k-\frac{\gamma}{2}
$$
for all $t$ sufficiently large. This indeed holds as the first term in the left hand side tends to $e^k$ as $t \rightarrow \infty$ while the two other terms on the left hand side tend to
$1$ as $t \rightarrow \infty$. \qed

It is easy to see that what really matters in the proof of Theorem \ref{t:k-reg} is not the $k$-regularity itself, but rather that $H$ has some bounded maximum degree and
that $e(H)=nk$ and also $plus(H)-minus(H)=c(H)+3f(H)-i(H)-g(H) = kn$. As for any even orientation (and hence any Eulerian orientation) we have $plus(H)-minus(H)=e(H)$ we obtain
the following corollary.
\begin{corollary}\label{c:final}
For any positive integer $K$ and any $\gamma > 0$, there exists $N=N(K,\gamma)$ such that if $H$ is any Eulerian orientation with $n > N$ vertices, maximum degree at most $K$, and $kn$ edges,
then $R(H) \ge (e^k-\gamma)n!/2^{nk}$.  \qed
\end{corollary}

\section{Concluding remarks}

{\bf Proof of Corollary \ref{c:main}.}
We slightly modify the definition of an adjusted $t$-decomposition to apply to the case where $n$ is even.
So, assume the vertices of $K_n$ are labeled with $[n]$ where $n$ is even. Take an adjusted $t$-decomposition ${\cal D}^*$ on $[n-1]$ (note: we still assume that $t$ is odd).
Now, to extend it to a decomposition of $K_n$ just add the following $n/2$ elements $P_1,\ldots,P_{n/2}$  where $P_i$ consists of the two edges $\{2i-1,n\}$ and $\{2i,n\}$
for $i=1,\ldots,n/2-1$ while $P_{n/2}$ is the single edge $\{n-1,n\}$. This decomposition ${\cal D}$ of $K_n$ is now defined as an adjusted $t$-decomposition.
Notice that each element of ${\cal D}$ that is not in ${\cal D}^*$ is either a $K_{1,2}$ or a single edge.
Now, in the proof of Theorem \ref{t:special} where we randomly orient each element of ${\cal D}$, do as before for the elements of ${\cal D}^*$, while for each $P_i$
for $i=1,\ldots,n/2-1$ orient it as a directed path of length $2$ in one of the two possible directions. Also, orient the single edge $P_{n/2}$ in one of the two possible
direction. Notice that the obtained random orientation $T$ induces a regular tournament on $[n-1]$ and vertex $n$ has in-degree either $n/2$ or $n/2-1$, so $T$
is a balanced tournament. Now, all the rest of the proof of Theorem \ref{t:special} carries over exactly in the the same way.
\qed

Recall that Wormald conjectured that $T(C_n)=\Theta(n!/2^n)$ and in fact that $T(C_n) \le 2.856n!/2^n$ for all sufficiently large $n$.
It would thus be natural to conjecture that $T(H)=\Theta(n!/2^{e(H)})$ for, say, all $(\epsilon,k)$-consistent orientations with $n$ vertices.

Perhaps less ambitious, but possibly difficult, is to improve the constant $e^k$ in Theorem \ref{t:k-reg} and in Corollary \ref{c:final}.
Recall that Wormald proved that for the Hamilton cycle (which applies the case $k=1$) the constant can be improved from $e$ to $2.855$.
However, it seems difficult to try to generalize Wormald's argument to all, say, $k$-regular orientations as such orientations may have highly involved structure
(they can contain, for example, arbitrary sub-tournaments on $\Theta(k)$ vertices) so it does not seem plausible to reduce the problem to a computer verification check of small constant
size tournaments an in the proof in \cite{wormald-p}.

Finally, recall that the results of \cite{alon-1990,FK-2005} prove that $T(C_n) = O(n^c n!/2^n)$ for a relatively small constant $c$.
It would be nice to prove an analogous result for say, $k$-regular orientations, stating $T(H) = O(n^c n!/2^{nk})$ for some constant $c=c(k)$.

\bibliographystyle{plain}

\bibliography{paper}

\end{document}